\documentclass[12pt,reqno]{amsart}

\usepackage[arrow,matrix,curve]{xy}

\usepackage{graphicx}         

\usepackage{amssymb, latexsym, amsmath, amscd, array,
}

\usepackage{float}           

\usepackage{hyperref}         
\usepackage[all]{hypcap}      

\newtheorem{theorem}{Theorem}[section]

\theoremstyle{definition}

\numberwithin{equation}{section}

\newcommand\N {{\mathbb N}}

\newcommand\R {{\mathbb R}}
\newcommand\Q {{\mathbb Q}}

\newcommand\st{{\rm st}}

\newcommand\RRR{\mbox{I\!I\!R}}
\newcommand\NNN{\mbox{I\!I\!N}}

\newcommand\Los{{\L}o{\'s}}

\DeclareMathOperator{\adequal}{\,{}_{\ulcorner\!\urcorner}\,}

\begin{document}

\thispagestyle{empty}

\title[Ten misconceptions and their debunking] {Ten misconceptions
from the history of analysis and their debunking}

\author[B\l aszczyk]{Piotr B\l aszczyk$^0$}

\author[Katz]{Mikhail G. Katz}

\author[Sherry]{David Sherry}

\address{Institute of Mathematics, Pedagogical University of Cracow,
Poland} \email{pb@ap.krakow.pl}

\address{Department of Mathematics, Bar Ilan University, Ramat Gan
52900 Israel} \email{katzmik@macs.biu.ac.il}

\address{Department of Philosophy, Northern Arizona University,
Flagstaff, AZ 86011, USA} \email{david.sherry@nau.edu}

\footnotetext{Supported by Polish Ministry of Science and Higher
Education grant N N101 287639}

\subjclass[2000]{%
01A85;            
Secondary
26E35,            
03A05,            
97A20,            
97C30             
}

\keywords{Abraham Robinson, adequality, Archimedean continuum,
Bernoullian continuum, Cantor, Cauchy, cognitive bias, completeness,
constructivism, continuity, continuum, du Bois-Reymond, epsilontics,
Felix Klein, Fermat-Robinson standard part, infinitesimal,
Leibniz-{\L}o{\'s} transfer principle, limit, mathematical rigor,
nominalism, non-Archimedean, Simon Stevin, Stolz, sum theorem,
Weierstrass}

\begin{abstract}
The widespread idea that infinitesimals were ``eliminated'' by the
``great triumvirate'' of Cantor, Dedekind, and Weierstrass, is refuted
by an uninterrupted chain of work on infinitesimal-enriched number
systems.  The elimination claim is an oversimplification created by
triumvirate followers, who tend to view the history of analysis as a
pre-ordained march toward the radiant future of Weierstrassian
epsilontics.  In the present text, we document distortions of the
history of analysis stemming from the triumvirate ideology of
ontological minimalism, which identified the continuum with a single
number system.  Such anachronistic distortions characterize the
received interpretation of Stevin, Leibniz, d'Alembert, Cauchy, and
others.
\end{abstract}

\maketitle

\tableofcontents

\section{Introduction}

Here are some common claims.  The founders of infinitesimal calculus
were working in a vacuum caused by an absence of a satisfactory number
system.  The incoherence of infinitesimals was effectively criticized
by Berkeley as so much hazy metaphysical mysticism.  D'Alembert's
visionary anticipation of the rigorisation of analysis was ahead of
his time.  Cauchy took first steps toward replacing infinitesimals by
rigor and epsilontics, in particular giving a modern definition of
continuity.  Cauchy's false 1821 version of the ``sum theorem" was
corrected by him in 1853 by adding the hypothesis of uniform
convergence.  Weierstrass finally rigorized analysis and thereby
eliminated infinitesimals from mathematics.  Dedekind discovered ``the
essence of continuity", which is captured by his cuts.  One of the
spectacular successes of the rigorous analysis was the mathematical
justification of Dirac's ``delta functions''.  Robinson developed a
new theory of infinitesimals in the 1960s, but his approach has little
to do with historical infinitesimals.  Lakatos pursued an ideological
agenda of Kuhnian relativism and fallibilism, inapplicable to
mathematics.

Each of the above ten claims is in error, as we argue in the next ten
sections (cf.~Crowe \cite{Cr}).

The historical fact of the dominance of the view of analysis as being
based on the real numbers to the exclusion of infinitesimals, is
beyond dispute.  One could ask oneself why this historical fact is so;
some authors have criticized mathematicians for adhering to an
approach that others consider less appropriate.  In the present text,
we will \emph{not} be concerned with either of these issues.  Rather,
we will be concerned with another issue, namely, why is it that
traditional historical scholarship has been inadequate in indicating
that alternative views have been around.  We will also be concerned
with documenting instances of tendentious interpretation and the
attendant distortion in traditional evaluation of key figures from
mathematical history.

Felix Klein clearly acknowledged the existence of a parallel,
infinitesimal approach to foundations.  Having outlined the
developments in real analysis associated with Weierstrass and his
followers, Klein pointed out in~1908 that
\begin{quote}
The scientific mathematics of today is built upon the series of
developments which we have been outlining.  But an essentially
different conception of infinitesimal calculus has been running
parallel with this [conception] through the centuries (Klein
\cite[p.~214]{Kl}).
\end{quote}
Such a different conception, according to Klein, ``harks back to old
metaphysical speculations concerning the structure of the continuum
according to which this was made up of [\ldots] infinitely small
parts'' (ibid.).  The pair of parallel conceptions of analysis are
illustrated in Figure~\ref{31}.

\begin{figure}
\[
\xymatrix@C=95pt{{} \ar@{-}[rr] \ar@{-}@<-0.5pt>[rr]
\ar@{-}@<0.5pt>[rr] & {} \ar@{->}[d]^{\hbox{st}} & \hbox{\quad
B-continuum} \\ {} \ar@{-}[rr] & {} & \hbox{\quad A-continuum} }
\]
\caption{\textsf{Pair of parallel conceptions of the continuum (the
thickness of the top line is merely conventional).  The
``thick-to-thin'' vertical arrow ``st'' represents taking the standard
part (see Appendix~\ref{rival2} for details).}}
\label{31}
\end{figure}

A comprehensive re-evaluation of the history of infinitesimal calculus
and analysis was initiated by Katz \& Katz in \cite{KK11a},
\cite{KK11b}, and~\cite{KK11c}.  Briefly, a philosophical disposition
characterized by a preference for a sparse ontology has dominated the
historiography of mathematics for the past 140 years, resulting in a
systematic distortion in the interpretation of the historical
development of mathematics from Stevin (see \cite{KK11c}) to Cauchy
(see \cite{KK11b} and Borovik \& Katz \cite{BK}) and beyond.  Taken to
its logical conclusion, such distortion can assume comical
proportions.  Thus, Newton's eventual successor in the Lucasian chair
of mathematics, Stephen Hawking, comments that Cauchy
\begin{quote}
was particularly concerned to banish infinitesimals (Hawking
\cite[p.~639]{Ha}),
\end{quote}
yet {\em on the very same page\/} 639, Hawking quotes Cauchy's {\em
infinitesimal\/} definition of continuity in the following terms:
\begin{quote}
the function~$f(x)$ remains continuous with respect to~$x$ between the
given bounds, if, between these bounds, an infinitely small increment
in the variable always produces an infinitely small increment in the
function itself (ibid).
\end{quote}
Did Cauchy banish infinitesimals, or did he exploit them to define a
seminal new notion of continuity?  Similarly, historian J.~Gray lists
{\em continuity\/} among concepts Cauchy allegedly defined
\begin{quote}
using careful, if not altogether unambiguous, {\bf limiting} arguments
(Gray \cite[p.~62]{Gray08}) [emphasis added--authors],
\end{quote}
whereas in reality {\em limits\/} appear in Cauchy's definition only
in the sense of the {\em endpoints\/} of the domain of definition (see
\cite{KK11a}, \cite{KK11b} for a more detailed discussion).
Commenting on `Whig' re-writing of mathematical history,%
\footnote{\label{GG1}Related comments by Grattan-Guinness may be found
in the main text at footnote~\ref{muddle}.}
P.~Mancosu observed that
\begin{quote}
the literature on infinity is replete with such `Whig' history.
Praise and blame are passed depending on whether or not an author
might have anticipated Cantor and naturally this leads to a completely
anachronistic reading of many of the medieval and later contributions
\cite[p.~626]{Ma09}.
\end{quote}
The anachronistic idea of the history of analysis as a relentless
march toward the yawning heights of epsilontics is, similarly, our
target in the present text.  We outline some of the main distortions,
focusing on the philosophical bias which led to them.  The outline
serves as a program for further investigation.

\section{Were the founders of calculus working in a numerical vacuum?}
\label{stev}


Were the founders of infinitesimal calculus working in a vacuum caused
by an absence of a satisfactory number system?

\subsection{Stevin, \emph{La Disme}, and \emph{Arithmetique}}

A century before Newton and Leibniz, Simon Stevin (Stevinus) sought to
break with an ancient Greek heritage of working exclusively with
relations among natural numbers,%
\footnote{The Greeks counted as numbers only~$2, 3, 4, \ldots$; thus,
$1$ was not a number, nor are the fractions: ratios were relations,
not numbers.  Consequently, Stevin had to spend time arguing that the
unit
%
%
($\mu o\nu\acute\alpha\varsigma$) {\em was\/} a number.}
and developed an approach capable of representing both ``discrete''
number
%
%
(\,'\!\!\!$\alpha\rho\iota\theta\mu \acute o\varsigma$)%
\footnote{Euclid \cite[Book VII, def. 2]{Euc}.}
composed of units 
%
%
($\mu o\nu\acute\alpha\delta\omega\nu$) and continuous magnitude
%
%
($\mu\acute\varepsilon\gamma\varepsilon\theta o\varsigma$) of
geometric origin.%
\footnote{Euclid \cite[Book V]{Euc}.  See also Aristotle's
\emph{Categories}, 6.4b, 20-23: ``Quantity is either discrete or
continuous. [\ldots] Instances of discrete quantities are number and
speech; of continuous, lines, surfaces, solids, and besides these,
time and place".}
According to van der Waerden, Stevin's
\begin{quote}
general notion of a real number was accepted, tacitly or explicitly,
by all later scientists \cite[p.~69]{van}.
\end{quote}
D.~Fearnley-Sander wrote that
\begin{quote}
the modern concept of real number [...] was essentially achieved by
Simon Stevin, around 1600, and was thoroughly assimilated into
mathematics in the following two centuries \cite[p.~809]{Fea}.
\end{quote}
D.~Fowler points out that
\begin{quote}
Stevin [...] was a thorough-going arithmetizer: he published, in 1585,
the first popularization of decimal fractions in the West [...]; in
1594, he described an algorithm for finding the decimal expansion of
the root of any polynomial, the same algorithm we find later in
Cauchy's proof of the Intermediate Value Theorem \cite[p.~733]{Fo}.
\end{quote}

Fowler \cite{Fo} emphasizes that important foundational work was yet
to be done by Dedekind, who proved that the field operations and other
arithmetic operations extend from~$\Q$ to~$\R$ (see
Section~\ref{dede}).%
\footnote{Namely, there is no easy way from representation of reals by
decimals, to the {\em field\/} of reals, just as there is no easy way
from continuous fractions, another well-known representation of reals,
to operations on such fractions.}
Meanwhile, Stevin's decimals stimulated the emergence of power series
(see below) and other developments.  We will discuss Stevin's
contribution to the Intermediate Value Theorem in Subsection~\ref{IVT}
below.

In 1585, Stevin defined decimals in \emph{La Disme} as follows:
\begin{quote}
Decimal numbers are a kind of arithmetic based on the idea of the
progression of tens, making use of the Arabic numerals in which any
number may be written and by which all computations that are met in
business may be performed by integers alone without the aid of
fraction. (\emph{La Disme}, On Decimal Fractions, tr. by V. Sanford,
in Smith \cite[p.~23]{Smith}).
\end{quote}
  
By numbers ``met in business" Stevin meant finite decimals,%
\footnote{\label{ad1}But see Stevin's comments on extending the
process \emph{ad infinitum} in main text at footnote~\ref{ad2}.}
and by ``computations'' he meant addition, subtraction,
multiplications, division and extraction of square roots on finite
decimals.%
\footnote{An algorithmic approach to such operations on infinite
decimals was developed by Hoborski \cite{Ho}, and later in a very
different way by Faltin et al.~\cite{Fa}.}
Stevin argued that numbers, like the more familiar continuous
magnitudes, can be divided indefinitely, and used a water metaphor to
illustrate such an analogy:
\begin{quote}
As to a continuous body of water corresponds a continuous wetness, so
to a continuous magnitude corresponds a continuous number.  Likewise,
as the continuous body of water is subject to the same division and
separation as the water, so the continuous number is subject to the
same division and separation as its magnitude, in such a way that
these two quantities cannot be distinguished by continuity and
discontinuity (Stevin, 1585, see \cite[p.~3]{St85}; quoted in A.~Malet
\cite{Mal06}).%
\footnote{See also Naets \cite{Nae} for an illuminating discussion of
Stevin.}
\end{quote}
Stevin argued for equal rights in his system for rational and
irrational numbers.  He was critical of the complications in Euclid
\cite[Book X]{Euc}, and was able to show that adopting the arithmetic
as a way of dealing with those theorems made many of them easy to
understand and easy to prove.  In his {\em Arithmetique\/}
\cite{St85}, Stevin proposed to represent all numbers systematically
in decimal notation.  P.~Ehrlich notes that Stevin's
\begin{quote}
viewpoint soon led to, and was implicit in, the analytic geometry of
Ren\'e Descartes (1596-1650), and was made explicit by John Wallis
(1616-1703) and Isaac Newton (1643-1727) in their arithmetizations
thereof (Ehrlich \cite[p.~494]{Eh05}).  
\end{quote}

\subsection{Decimals from Stevin to Newton and Euler}

Stevin's text \emph{La Disme} on decimal notation was translated into
English in 1608 by Robert Norton (cf.~Cajori \cite[p.~314]{Caj}).  The
translation contains the first occurrence of the word ``decimal'' in
English; the word will be employed by Newton in a crucial passage 63
years later.%
\footnote{\label{new1}See footnote~\ref{new2}.}
Wallis recognized the importance of unending decimal expressions in
the following terms:
\begin{quote}
Now though the Proportion cannot be accurately expressed in absolute
Numbers: Yet by continued Approximation it may; so as to approach
nearer to it, than any difference assignable (Wallis's Algebra,
p. 317, cited in Crossley \cite{Cr87}).
\end{quote}
Similarly, Newton exploits a power series expansion to calculate
detailed decimal approximations to~$log(1+x)$ for~$x=\pm 0.1, \pm 0.2,
\dots$ (Newton \cite{New??}).%
\footnote{Newton scholar N.~Guicciardini kindly provided a jpg of a
page from Newton's manuscript containing such calculations, and
commented as follows: ``It was probably written in Autumn 1665 (see
Mathematical papers 1, p.~134).  Whiteside's dating is sometimes too
precise, but in any case it is a manuscript that was certainly written
in the mid 1660s when Newton began annotating Wallis's
\emph{Arithmetica Infinitorum}'' \cite{Gu}.  The page from Newton's
manuscript can be viewed at
\url{http://u.math.biu.ac.il/~katzmik/newton.html }}
By the time of Newton's \emph{annus mirabilis}, the idea of unending
decimal representation was well established.  Historian V.~Katz calls
attention to ``Newton's analogy of power series to infinite decimal
expansions of numbers'' (V.~Katz \cite[p.~245]{KV}).  Newton expressed
such an analogy in the following passage:
\begin{quote}
Since the operations of computing in numbers and with variables are
closely similar \ldots I am amazed that it has occurred to no one (if
you except N. Mercator with his quadrature of the hyperbola) to fit
the doctrine recently established for decimal numbers in similar
fashion to variables, especially since the way is then open to more
striking consequences.  For since this doctrine in species has the
same relationship to Algebra that the doctrine in decimal numbers has
to common Arithmetic, its operations of Addition, Subtraction,
Multiplication, Division and Root-extraction may easily be learnt from
the latter's provided the reader be skilled in each, both Arithmetic
and Algebra, and appreciate the correspondence between decimal numbers
and algebraic terms continued to infinity \ldots And just as the
advantage of decimals consists in this, that when all fractions and
roots have been reduced to them they take on in a certain measure the
nature of integers, so it is the advantage of infinite
variable-sequences that classes of more complicated terms (such as
fractions whose denominators are complex quantities, the roots of
complex quantities and the roots of affected equations) may be reduced
to the class of simple ones: that is, to infinite series of fractions
having simple numerators and denominators and without the all but
insuperable encumbrances which beset the others (Newton 1671,
\cite{New71}).%
\footnote{We are grateful to V.~Katz for signaling this passage.}
\end{quote}
In this remarkable passage dating from 1671, Newton explicitly names
infinite decimals as the source of inspiration for the new idea of
infinite series.%
\footnote{\label{new2}See footnote~\ref{new1}.}
The passage shows that Newton had an adequate number system for doing
calculus and real analysis two centuries before the triumvirate.%
\footnote{The expression ``the great triumvirate'' is used by Boyer
\cite[p.~298]{Boy} to describe Cantor, Dedekind, and Weierstrass.}

In 1742, John Marsh first used an abbreviated notation for repeating
decimals (Marsh \cite[p.~5]{Mar}, cf.~Cajori \cite[p.~335]{Caj}).
Euler exploits unending decimals in his {\em Elements of Algebra\/} in
1765, as when he sums an infinite series and concludes%
\footnote{\label{per1}This could be compared with Peirce's remarks.
Over a century ago, Charles Sanders Peirce wrote with reference
to~$1-\text{`}0.999\ldots$': ``although the difference, being
infinitesimal, is less than any number [one] can express[,] the
difference exists all the same, and sometimes takes a quite easily
intelligible form'' (Peirce \cite[p.~597]{Pe}; see also
S.~Levy~\cite[p.~130]{Levy}).  Levy mentions Peirce's proposal of an
alternative notation for ``equality up to an infinitesimal".  The
notation Peirce proposes is the usual equality sign with a dot over
it, like this:~``$\dot=$''.  See also main text at
footnote~\ref{per2}.}
that~$9.999\ldots=10$ (Euler \cite[p.~170]{Eul}).

\subsection{Stevin's cubic and the IVT}
\label{IVT}

In the context of his decimal representation, Stevin developed
numerical methods for finding roots of polynomials equations.  He
described an algorithm equivalent to finding zeros of polynomials (see
Crossley \cite[p.~96]{Cr87}).  This occurs in a corollary to problem
77 (more precisely, LXXVII) in (Stevin \cite[p.~353]{St1625}).  Here
Stevin describes such an argument in the context of finding a root of
the cubic equation (which he expresses as a proportion to conform to
the contemporary custom)
\[
x^3=300x+33915024.
\]
Here the whimsical coefficient seems to have been chosen to emphasize
the fact that the method is completely general; Stevin notes
furthermore that numerous addional examples can be given.  A textual
discussion of the method may be found in Struik~\cite[p.~476]{Ste}.

Centuries later, Cauchy would prove the Intermediate Value Theorem
(IVT) for a continuous function~$f$ on an interval~$I$ by a
divide-and-conquer algorithm.  Cauchy subdivided~$I$ into~$m$ equal
subintervals, and recursively picked a subinterval where the values of
$f$ have opposite signs at the endpoints (Cauchy \cite[Note III,
p.~462]{Ca21}).  To elaborate on Stevin's argument following
\cite[\S10, p.~475-476]{Ste}, note that what Stevin similarly
described a divide-and-conquer algorithm.  Stevin subdivides the
interval into {\em ten\/} equal parts, resulting in a gain of a new
decimal digit of the solution at every iteration of his procedure.
Stevin explicitly speaks of continuing the iteration \emph{ad
infinitum}:%
\footnote{\label{ad2}Cf. footnote~\ref{ad1}.}
\begin{quote}
Et procedant ainsi infiniment, l'on approche infiniment plus pres au
requis (Stevin \cite[p.~353, last 3 lines]{St1625}).
\end{quote}

Who needs {\em existence\/} proofs for the real numbers, when Stevin
gives a procedure seeking to produce an explicit decimal
representation of the solution?  The IVT for polynomials would
resurface in Lagrange before being generalized by Cauchy to the newly
introduced class of continuous functions.%
\footnote{\label{lag1}See further in footnote~\ref{lag2}.}

One frequently hears sentiments to the effect that the pre-triumvirate
mathematicians did not and could not have provided rigorous proofs,
since the real number system was not even built yet.  Such an attitude
is anachronistic.  It overemphasizes the significance of the
triumvirate project in an inappropriate context.  Stevin is concerned
with constructing an algorithm, whereas Cantor is concerned with
developing a foundational framework based upon the existence of the
{\em totality\/} of the real numbers, as well as their power sets,
etc.  The latter project involves a number of non-constructive
ingredients, including the axiom of infinity and the law of excluded
middle.  But none of it is needed for Stevin's procedure, because he
is not seeking to re-define ``number'' in terms of alternative
(supposedly less troublesome) mathematical objects.

Why do many historians and mathematicians of today emphasize the great
triumvirate's approach to proofs of the existence of real numbers, at
the expense, and almost to the exclusion, of Stevin's approach?  Can
this be understood in the context of the ideological foundational
battles raging at the end of 19th and beginning of 20th century?
These questions merit further scrutiny.

\section{Was Berkeley's criticism coherent?}


Was Berkeley's criticism of infinitesimals as so much hazy
metaphysical mysticism, either effective or coherent?

D.~Sherry \cite{She87} dissects Berkeley's criticism of infinitesimal
calculus into its metaphysical and logical components, as detailed
below.

\subsection{Logical criticism}
The {\em logical criticism\/} is the one about the disappearing~$dx$.
Here we have a ghost:~$dx\not=0$, but also a departed quantity:~$dx=0$
(in other words eating your cake:~$dx=0$ and having it,
too:~$dx\not=0$).

Thus, Berkeley's {\em logical criticism\/} of the calculus is that the
evanescent increment is first assumed to be non-zero to set up an
algebraic expression, and then {\em treated as zero\/} in {\em
discarding\/} the terms that contained that increment when the
increment is said to vanish.%
\footnote{Given Berkeley's fame among historians of mathematics for
allegedly spotting logical flaws in infinitesimal calculus, it is
startling to spot circular logic at the root of Berkeley's own
doctrine of the compensation of errors.  Indeed, Berkeley's new,
improved calculation of the derivative of $x^2$ in \emph{The Analyst}
\cite{Be} relies upon the determination of the tangent to a parabola
due to Apollonius of Perga \cite[Book I, Theorem 33]{Ap} (see Andersen
\cite{An11}).}

The fact is that Berkeley's logical criticism is easily answered
within a conceptual framework available to the founders of the
calculus.  Namely, the {\em rebuttal\/} of the logical criticism is
that the evanescent increment is not {\em treated as zero\/}, but,
rather, merely {\em discarded\/} through an application of Leibniz's
{\em law of homogeneity\/} (see Leibniz \cite[p.~380]{Le10b}) and Bos
\cite[p.~33]{Bos}), which would stipulate, for instance, that
\begin{equation}
\label{31b}
a+dx\adequal a.
\end{equation}
Here we chose the sign~$\adequal$ which was already used by Leibniz
where we would use an equality sign today (see McClenon
\cite[p.~371]{Mc23}).  The law is closely related to the earlier
notion of {\em adequality\/} found in Fermat.  Adequality is the
relation of being infinitely close, or being equal ``up to'' an
infinitesimal.  Fermat exploited adequality when he sought a maximum
of an expression by evaluating expression at~$A+E$ and at~$A$, and
forming the difference.  In modern notation this would appear
as~$f(A+E)-f(A)$ (note that Fermat did not use the function notation).
Huygens already interpreted the ``$E$'' occurring in this expression
in the method of adequality, as an infinitesimal.%
\footnote{Huygens explained Fermat's method of adequality in a
presentation at the {\em Acad\'emie des Sciences\/} in 1667.  Huygens
noted that ``$E''$ is an ``infinitely small quantity'' (see Huygens
\cite{Hu}).  See also Andr\'e Weil \cite[p.~1146]{We},
\cite[p.~28]{We84}.}

\begin{figure}
\includegraphics[height=2.1in]{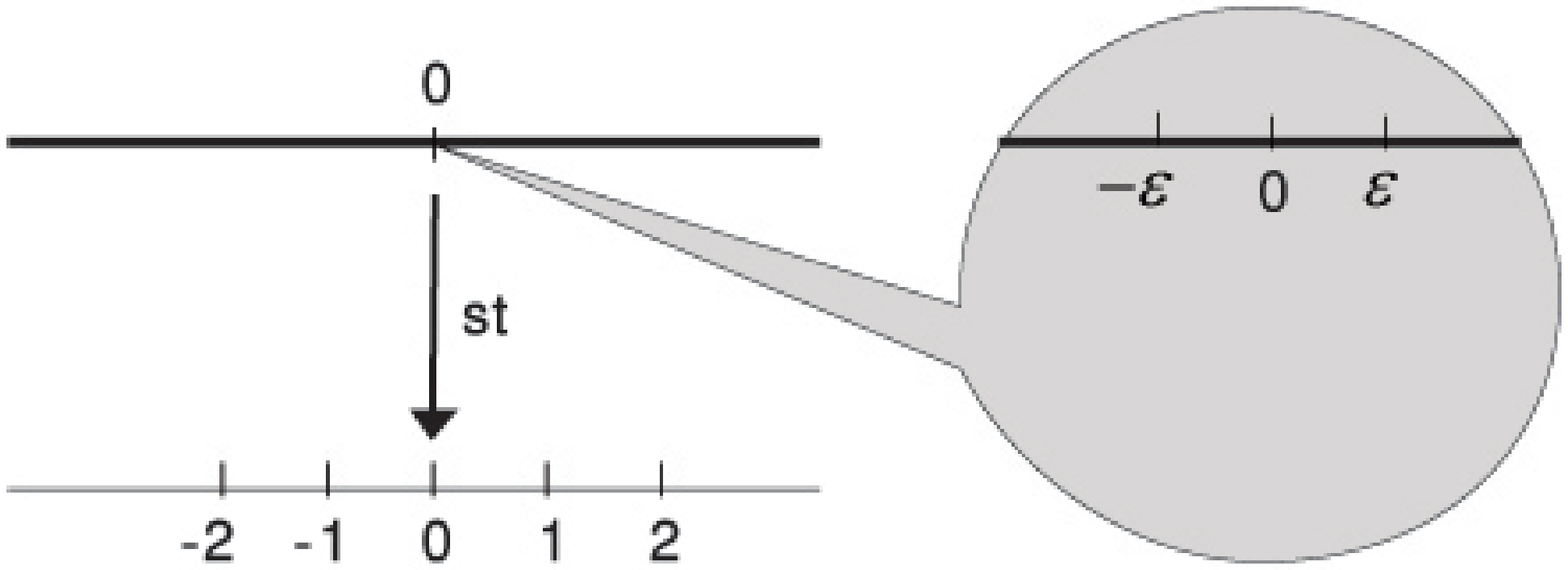}
\caption{\textsf{Zooming in on infinitesimal~$\varepsilon$ (see
Appendix~\ref{rival2} for details)}}
\label{tamar}
\end{figure}

Ultimately, the heuristic concepts of adequality (Fermat) and law of
homogeneity (Leibniz) were implemented in terms of the {\em standard
part\/} function (see Figure~\ref{tamar}).

In a passage typical of post-Weierstrassian scholarship, Kleiner and
Movshovitz-Hadar note that
\begin{quote}
Fermat's method was severely criticized by some of his contemporaries.
They objected to his introduction and subsequent suppression of the
mysterious~$E$.  Dividing by~$E$ meant regarding it as not zero.
Discarding~$E$ implied treating it as zero.  This is inadmissible,
{\bf they rightly claimed.}  In a somewhat different context, but {\bf
with equal justification}, \ldots{} Berkeley in the 18th century would
refer to such~$E$'s as `the ghosts of departed quantities' ''
\cite[p.~970]{KM} [emphasis added--authors].
\end{quote}

However, Fermat scholar P.~Str\o mholm already pointed out in 1968
that in Fermat's main method of adequality,
\begin{quote}
there was never [\ldots] any question of the variation~$E$ being put
equal to zero.  The words Fermat used to express the process of
suppressing terms containing~$E$ was {\em ``elido''\/}, {\em
``deleo''\/}, and {\em ``expungo''\/}, and in French {\em
``i'efface''\/} and {\em ``i'\^ote''\/}.  We can hardly believe that a
sane man wishing to express his meaning and searching for words, would
constantly hit upon such tortuous ways of imparting the simple fact
that the terms vanished because~$E$ was zero (Str\o mholm
\cite[p.~51]{Strom}).
\end{quote}

Thus, Fermat planted the seeds of the answer to the logical criticism
of the infinitesimal, a century before George Berkeley ever lifted up
his pen to write {\em The Analyst\/}.

The existence of two separate binary relations, one ``equality'' and
the other, ``equality up to an infinitesimal'', was already known to
Leibniz (see E.~Knobloch \cite[p.~63]{Kn} and Katz and Sherry
\cite{KS} for more details).

\subsection{Metaphysical criticism}
Berkeley's {\em metaphysical criticism\/} targets the absence of any
empirical referent for ``infinitesimal''.  The metaphysical criticism
has its roots in empiricist dogma that every meaningful expression or
symbol must correspond to an empirical entity.%
\footnote{The interplay of empiricism and nominalism in Berkeley's
thought is touched upon by D.~Sepkoski \cite[p.~50]{Se}.}
Ironically, Berkeley accepts many expressions lacking an empirical
referent, such as `force', `number', or `grace', on the grounds that
they have pragmatic value.  It is a glaring inconsistency on
Berkeley's part not to have accepted ``infinitesimal" on the same
grounds (see Sherry~\cite{She93}).

\begin{figure}
\includegraphics[height=2.2in]{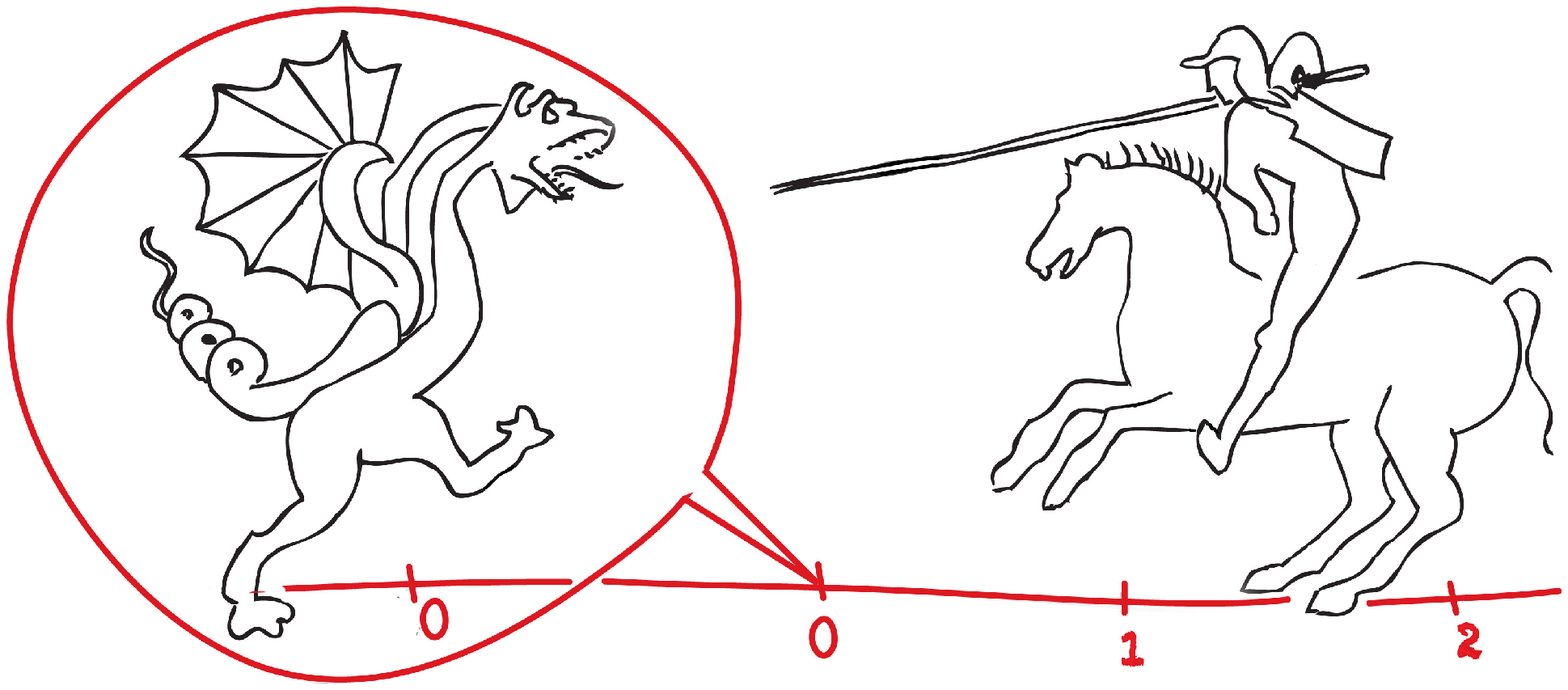}
\caption{\textsf{George's attempted slaying of the infinitesimal,
following E.~T.~Bell and P.~Uccello.  Uccello's creature is shown as
inhabiting an infinitesimal neighborhood of~$0$,
cf.~Figure~\ref{tamar}.}}
\label{uccello}
\end{figure}

It is even more ironic that over the centuries, mathematicians were
mainly unhappy with the logical aspect, but their criticisms mainly
targeted what they perceived as the metaphysical/mystical aspect.
Thus, Cantor attacked infinitesimals as being ``abominations" (see
Ehrlich \cite{Eh06}); R.~Courant described them as ``mystical", ``hazy
fog'', etc.  E.~T. Bell went as far as describing infinitesimals as
having been
\begin{itemize}
\item
{\em slain\/} \cite[p.~246]{Bel45},
\item
{\em scalped\/} \cite[p.~247]{Bel45}, and
\item
{\em disposed of\/} \cite[p.~290]{Bel45}
\end{itemize}
by the cleric of Cloyne (see Figure~\ref{uccello}).  Generally
speaking, one does not {\em slay\/} either scientific concepts or
scientific entities.  Bellicose language of this sort is a sign of
commitments that are both emotional and ideological.

\section{Were d'Alembert's anticipations ahead of his time?}
\label{four}


Were d'Alembert's mathematical anticipations ahead or behind his time?

One aspect of d'Alembert's vision as expressed in his article for the
{\em Encyclopedie\/} on {\em irrational\/} numbers, is that irrational
numbers \emph{do not exist}.  Here d'Alembert uses terms such as
``surds" which had already been rejected by Simon Stevin two centuries
earlier (see Section~\ref{stev}).  From this point of view, d'Alembert
is not a pioneer of the rigorisation of analysis in the 19th century,
but on the contrary, represents a throwback to the 16th century.
D'Alembert's attitude toward irrational numbers sheds light on the
errors in his proof of the fundamental theorem of algebra;%
\footnote{D'Alembert's error was first noticed by Gauss, who gave some
correct proofs, though S. Smale argues his first proof was incomplete
\cite{Sma}.  See Baltus \cite{Bal} for a detailed study of
d'Alembert's proof.}
indeed, in the anemic number system envisioned by d'Alembert, numerous
polynomials surely fail to have roots.

D'Alembert used the term ``charlatanerie" to describe infinitesimals
in his article {\em Diff\'erentiel\/} \cite{da}.  D'Alembert's
anti-infinitesimal vitriol is what endears him to triumvirate
scholars, for his allegedly visionary remarks concerning the
centrality of the limit concept fall short of what is already found in
Newton.%
\footnote{\label{pourciau}See Pourciau \cite{Pou} who argues that
Newton possessed a clear kinetic conception of limit (Sinaceur
\cite{Si} and Barany \cite{Ba11} argue that Cauchy's notion of limit
was kinetic, rather than a precursor of a Weierstrassian notion of
limit).  Pourciau cites Newton's lucid statement to the effect that
``Those ultimate ratios \ldots are not actually ratios of ultimate
quantities, but limits \ldots which they can approach so closely that
their difference is less than any given quantity\ldots'' (Newton, 1946
\cite[p.~39]{New46} and 1999 \cite[p.~442]{New99}).  The same point,
and the same passage from Newton, appeared a century earlier in
Russell \cite[item 316, p.~338-339]{Ru03}.}
He never went beyond a kinetic notion of limit, so as to obtain the
epsilontic version popularized in the 1870s.

D'Alembert was particularly bothered by the characterisation of
infinitesimals he found in ``the geometers".  He does not explain who
these geometers are, but the characterisation he is objecting to can
be already found in Leibniz and Newton.  Namely, the geometers used to
describe infinitesimals as what remains
\begin{quote}
``not before you pass to the limit, nor after, but at the very moment
of passage to the limit".
\end{quote}

In the context of a modern theory of infinitesimals such as the
hyperreals (see Appendix~\ref{rival2}), one could explain the matter
in the following terms.  We decompose the procedure of taking the
limit of, say, a sequence~$(r_n)$ into two stages:
\begin{enumerate}
\item[(i)] evaluating the sequence at an infinite hypernatural%
\footnote{\label{hyper1}See main text at footnote~\ref{hyper2}.}
~$H$, to
obtain the hyperreal~$r_H$; and
\item[(ii)] taking its standard part~$L=\text{st}(r_H)$.
\end{enumerate}
Thus~$r_H$ is adequal to~$L$, or~$r_H\adequal L$ (see
formula~\eqref{31b}), so that we have~$\lim\limits_{n\to\infty}
r_n=L$.  In this sense, the infinitesimals exist ``at the moment" of
taking the limit, namely {\em between\/} the stages (i) and (ii).

Felscher \cite{Fe} describes d'Alembert as ``one of the mathematicians
representing the heroic age of calculus'' \cite[p.~845]{Fe}.  Felscher
buttresses his claim by a lengthy quotation concerning the definition
of the limit concept, from the article {\em Limite\/} from the {\em
Encyclop\'edie ou Dictionnaire Raisonn\'e des Sciences, des Arts et
des M\'etiers\/}:
\begin{quote}
On dit qu'une grandeur est la limite d'une autre grandeur, quand la
seconde peut approcher de la premi\`ere plus pr\`es que d'une grandeur
donn\'ee, si petite qu'on la puisse supposer, sans pourtant que la
grandeur qui approche, puisse jamais surpasser la grandeur dont elle
approche; ensorte que la diff\'erence d'une pareille quantit\'e \`a sa
limite est absolument inassignable (Encyclop\'edie, volume 9 from
1765, page 542).
\end{quote}
One recognizes here a kinetic definition of limit already exploited by
Newton.%
\footnote{See footnote~\ref{pourciau} on Pourciau's analysis.}
Even if we do attribute visionary status to this passage as many
historians seem to, the fact remains that d'Alembert didn't write it.
Felscher overlooked the fact that the article {\em Limite\/} was
written by two authors.  In reality, the above passage defining the
concept of ``limit" (as well as the two propositions on limits) did
not originate with d'Alembert, but rather with the encyclopedist
Jean-Baptiste de La Chapelle.  De la Chapelle was recruited by
d'Alembert to write 270 articles for the {\em Encyclop\'edie\/}.  The
section of the article containing these items is signed (E) (at bottom
of first column of page 542), known to be de La Chapelle's
``signature'' in the {\em Encyclopedie\/}.  Felscher had already
committed a similar error of attributing de la Chapelle's work to
d'Alembert, in his 1979 work \cite{Fe79}.%
\footnote{We are grateful to D.~Spalt for the historical clarification
concerning the authorship of the {\em Limite\/} article in the {\em
Encyclopedie\/}.}
Note that Robinson \cite[p.~267]{Ro66} similarly misattributes this
passage to d'Alembert.

\section{Did Cauchy replace infinitesimals by rigor?}
\label{did}


Did Cauchy take first steps toward replacing infinitesimals by rigor,
and did he give an epsilontic definition of continuity?

A claim to the effect that Cauchy was a fore-runner of the
epsilontisation of analysis is routinely recycled in history books and
textbooks.  To put such claims in historical perspective, it may be
useful to recall Grattan-Guinness's articulation of a historical
reconstruction project in the name of H.~Freudenthal~\cite{Fr71b}, in
the following terms:
\begin{quote}
it is mere feedback-style ahistory to read Cauchy (and contemporaries
such as Bernard Bolzano) as if they had read Weierstrass already.  On
the contrary, their own pre-Weierstrassian muddles%
\footnote{\label{muddle}Grattan-Guinness's term ``muddle'' refers to
an irreducible ambiguity of historical mathematics such as Cauchy's
sum theorem of 1821.  See footnote~\ref{GG1} for a related comment by
Mancosu.}
need historical reconstruction \cite[p.~176]{Grat04}.
\end{quote}

\subsection{Cauchy's definition of continuity}
\label{lutz}

It is often claimed that Cauchy gave an allegedly ``modern'', meaning
epsilon-delta, definition of continuity.  Such claims are
anachronistic.  In reality, Cauchy's definition is an infinitesimal
one.  His definition of the continuity of~$y=f(x)$ takes the following
form: an infinitesimal~$x$-increment gives rise to an
infinitesimal~$y$-increment (see \cite[p.~34]{Ca21}).  The widespead
misconception that Cauchy gave an epsilontic definition of continuity
is analyzed in detail in \cite{KK11a}.

Cauchy's primary notion is that of a {\em variable quantity\/}.  The
meaning he attached to the latter term in his {\em Cours d'Analyse\/}
in 1821 is generally agreed to have been a sequence of discrete
values.  He defines infinitesimals in terms of variable quantities, by
specifying that a variable quantity tending to zero {\em becomes\/} an
infinitesimal.  He similarly defines limits in terms of variable
quantities in the following terms:
\begin{quote}
lorsque les valeurs successivement attribu\'ees \`a une m\^eme
variable s'approche ind\'efiniment d'une valeur fixe de mani\`ere \`a
finir par en diff\'erer aussi peu que l'on voudra cette derni\`ere est
appel\'ee limite de toutes les autres (Cauchy \cite[p.~4]{Ca21}).
\end{quote}

Cauchy's definition is patently a kinetic, not an epsilontic,
definition of limit, similar to Newton's.%
\footnote{See footnote~\ref{pourciau}.}
While epsilontic-style formulations do occasionally appear in Cauchy
(though without Bolzano's proper attention to the order of the
quantifiers), they are not presented as definitions but rather as
consequences, and at any rate never appear in the context of the
property of the continuity of functions.

Thus, Grabiner's famous essay {\em Who gave you the epsilon?  Cauchy
and the origins of rigorous calculus\/} cites a form of an
epsilon-delta quantifier technique used by Cauchy in a proof:
\begin{quote}
Let~$\epsilon, \delta$ be two very small numbers; the first is chosen
so that for all numerical [i.e., absolute] values of~$h$ less than
$\delta$, and for any value of~$x$ included [in the interval of
definition], the ratio~$(f(x + h) - f(x))/h$ will always be greater
than~$f'(x) - \epsilon$ and less than~$f'(x) + \epsilon$ (Grabiner
\cite[p.~185]{Grab} citing Cauchy).
\end{quote}
Grabiner describes such an epsilon-delta technique as ``the algebra of
inequalities''.  The thrust of her argument is that Cauchy sought to
establish a foundation for analysis based on the algebra of
inequalities.  Is this borne out by the evidence she presents?  Let us
consider Grabiner's evidence:
\begin{quote}
Cauchy gave essentially the modern definition of continuous function,
saying that the function~$f(x)$ is continuous on a given interval if
for each~$x$ in that interval ``the numerical [i.e., absolute] value
of the difference~$f(x + \alpha) - f(x)$ decreases indefinitely with
$\alpha$ \cite[p.~190]{Grab}.
\end{quote}
Is this ``essentially the modern definition of continuity'', as
Grabiner claims?  Hardly so.  Cauchy's definition is a blend of a
kinetic (rather than epsilontic) and an infinitesimal approach.
Grabiner fails to mention three essential items:
\begin{itemize}
\item
Cauchy prefaces the definition she cited, by describing his~$\alpha$
as an {\em infinitely small increment\/}: 
\begin{quote}
Si, en partant d'une valeur de $x$ \ldots on attribue \`a la
variable~$x$ un \emph{accroissement infiniment petit}~$\alpha$
\ldots'' (Cauchy \cite[p.~34]{Ca21}) [emphasis added--the authors];
\end{quote}
\item
Cauchy follows this definition by another, {\em italicized\/},
definition, where both~$\alpha$ and the difference~$f(x+\alpha)-f(x)$
are described as being {\em infinitesimal\/}: if the former is
infinitesimal, then so is the latter;
\item
Infinitesimals provide a method for calculating limits, whereas
epsilon, delta methods require the answer in advance (see Madison and
Stroyan \cite[p.~497]{MS}).
\end{itemize}

The advantage of infinitesimal definitions, such as those found in
Cauchy, is their covariant nature (cf.~Lutz et al. \cite{LA}).
Whereas in the epsilontic approach one needs to work one's way
backwards from the value of the limit, in the infinitesimal approach
one can proceed from the original expression, simplify it, and
eventually arrive at the value of the limit.  This indicates that the
two approaches work in opposite directions.  The infinitesimal
calculation goes with the natural flow of our reasoning, whereas the
epsilontic one goes in the opposite direction.  Notice, for example,
that delta corresponds to the independent variable even though the
value of delta depends on our choice of epsilon, which corresponds to
the \emph{dependent} variable.  The infinitesimal calculation, in
contrast, begins with the the independent variable and computes from
that the value of the dependent variable.

\subsection{Cauchy's intermediate value theorem}

Did Cauchy exploit epsilon-delta techniques in building foundations
for analysis?  Let us examine Grabiner's evidence.  She claims that,
in Cauchy's proof of the intermediate value theorem (IVT), 
\begin{quote}
we have the algebra of inequalities providing a technique which Cauchy
transformed from a tool of approximation to a tool of rigor (Grabiner
\cite[p.~191]{Grab}).
\end{quote}
Yet Grabiner's treatment of Cauchy's proof of the IVT in
\cite[p.~190]{Grab} page offers no evidence that Cauchy employed an
epsilon-delta technique.%
\footnote{\label{lag2}Grabiner further attributes to Lagrange the
polynomial case of Cauchy's divide-and-conquer argument in the proof
of the IVT, whereas we saw in Subsection~\ref{IVT} that Stevin did
this two centuries before Lagrange (see main text at
footnote~\ref{lag1}).}

An examination of Cauchy's proof in his \emph{Note III} reveals that,
on the contrary, it is closely tied to Cauchy's infinitesimal
definition of continuity.  Thus, Cauchy constructs an increasing
sequence and a decreasing sequence, denoted respectively $x_n$ and
$X^{(n)}$ (Cauchy \cite[p.~462]{Ca21}) with a common limit $a$, such
that $f$ has opposite sign at the corresponding pairs of points.
Cauchy concludes that the values of $f$ at the respective sequences
converge to a common limit~$f(a)$.  Being both nonpositive and
nonnegative, the common limit $f(a)$ must vanish.

Koetsier \cite{Ko} speculates that Cauchy may have hit upon his
concept of continuity by analyzing his proof of the IVT (perhaps in
the case of polynomials).  The evidence is compelling: even though
Cauchy does not mention infinitesimals in his \emph{Note III}, $(x_n)$
and $\left(X^{(n)}\right)$ are recognizably variable quantities
differing by an infinitesimal from the constant quantity~$a$.  By
Cauchy's definition of continuity, $(f(x_n))$ and
$\left(f\left(X^{(n)}\right)\right)$ must similarly differ from $f(a)$
by an infinitesimal.  Contrary to Grabiner's claim, a close
examination of Cauchy's proof of the IVT reveals no trace of
epsilon-delta.  Following Koetsier's hypothesis, it is reasonable to
place it, rather, in the infinitesimal strand of the development of
analysis, rather than the epsilontic strand.

After constructing the lower and upper sequences, Cauchy does write
that the values of the latter ``finiront par differer de ces premiers
valeurs aussi peu que l'on voudra''.  That may sound a little bit
epsilon/delta.  Meanwhile, Leibniz uses language similar to Cauchy's:
\begin{quote}
Whenever it is said that a certain infinite series of numbers has a
sum, I am of the opinion that all that is being said is that any
finite series with the same rule has a sum, and that the error always
diminishes as the series increases, so that it becomes as small as we
would like [``ut fiat tam parvus quam velimus"] (Leibniz
\cite[p.~99]{Le01}).
\end{quote}
Cauchy used epsilontics if and only if Leibniz did, over a century
before him.

\subsection{Cauchy's influence}

The exaggerated claims of a Cauchy provenance for epsilontics found in
triumvirate literature go hand-in-hand with neglect of his visionary
role in the development of infinitesimals at the end of the 19th
century.  In 1902, E.~Borel \cite[p.~35-36]{Bo02} elaborated on Paul
du Bois-Reymond's theory of rates of growth, and outlined a general
``theory of increase'' of functions, as a way of implementing an
infinitesimal-enriched continuum.  In this text, Borel specifically
traced the lineage of such ideas to a 1829 text of Cauchy's
\cite{Ca29} on the rates of growth of functions (see Fisher
\cite[p.~144]{Fi} for details).  In 1966, A.~Robinson pointed out that
\begin{quote}
Following Cauchy's idea that an infinitely small or infinitely large
quantity is associated with the behavior of a function~$f(x)$, as~$x$
tends to a finite value or to infinity, du Bois-Reymond produced an
elaborate theory of orders of magnitude for the asymptotic behavior of
functions \dots Stolz tried to develop also a theory of arithmetical
operations for such entities \cite[p.~277-278]{Ro66}.
\end{quote}
Robinson traces the chain of influences further, in the following
terms:
\begin{quote}
It seems likely that Skolem's idea to represent infinitely large
natural numbers by number-theoretic functions which tend to infinity
(Skolem [1934]),%
\footnote{\label{1966a}The reference is to Skolem's 1934 work
\cite{Sk}.  The evolution of modern infinitesimals is traced in more
detail in Table~\ref{heuristic} and in Borovik et al.~\cite{BK}.}
also is related to the earlier ideas of Cauchy and du Bois-Reymond
\cite[p.~278]{Ro66}.
\end{quote}
One of Cantor's \emph{b\^etes noires} was the neo-Kantian philosopher
Hermann Cohen (1842--1918) (see also Mormann \cite{Mo}), whose
fascination with infinitesimals elicited fierce criticism by both
Cantor and B.~Russell.  Yet at the end of the day, A.~Fraenkel (of
Zermelo--Fraenkel fame) wrote:
\begin{quote}
my former student Abraham Robinson had succeeded in saving the honour
of infinitesimals - although in quite a different way than Cohen and
his school had imagined (Fraenkel \cite[p.~107]{Fra67}).
\end{quote}

\section{Was Cauchy's 1821 ``sum theorem" false?}
\label{sum}


Was Cauchy's 1821 ``sum theorem" false, and what did he add in 1853?

As discussed in Section~\ref{did}, Cauchy's definition of continuity
is explicitly stated in terms of infinitesimals: ``an
infinitesimal~$x$-increment gives rise to an infinitesimal
$y$-increment''.  Boyer \cite[p.~282]{Boy} declares that Cauchy's 1821
definition is ``to be interpreted" in the framework of the usual
``limits", at a point of an Archimedean continuum.  Traditional
historians typically follow Boyer's lead.

But when it comes to Cauchy's 1853 modification of the hypothesis of
the ``sum theorem"%
\footnote{The assertion of the theorem is the continuity of the sum of
a {\em convergent\/} series of continuous functions, with the
italicized term requiring clarification.  Modern versions of the
theorem require a hypothesis {\em uniform\/} convergence.  The nature
of the hypothesis Cauchy himself had in mind is hotly disputed (see
Borovik \& Katz \cite{BK} as well as \cite{KK11b}).}
in (Cauchy \cite{Ca53}), some historians declare that it is to be
interpreted as adding the hypothesis of ``uniform convergence" (see
e.g., L\"utzen \cite[p.~183-184]{Lut03}).
%
%
Are Boyer and L\"utzen compatible?

Note that an epsilontic definition (in the context of an Archimedean
continuum) of the uniform convergence of a sequence~$\langle f_n :
n\in\N \rangle$ to~$f$ necessarily involves a {\em pair\/} of
variables~$x,n$ (where~$x$ ranges through the domain of~$f$ and~$n$
ranges through~$\N$), rather than a single variable: we need a formula
of the sort
\begin{equation}
\label{32}
\forall n\in\N \; \forall x \;\left( n>n_0 \implies |f_n(x)-f(x)| <
\epsilon \right)
\end{equation}
(prefaced by the usual clause ``\,$(\forall \epsilon>0) \; (\exists
n_0\in\N)$\,'').  Now Cauchy's 1853 modification of the hypothesis is
stated in terms of a {\em single\/} variable~$x$, rather than a pair
of variables~$x,n$.  Namely, Cauchy specified that the condition of
convergence should hold ``always".  The meaning of the term ``always''
becomes apparent only in the course of the proof, when Cauchy gives an
explicit example of evaluating at an infinitesimal generated by the
sequence~$x=1/n$.  Thus the term ``always'' involves adding extra
values of~$x$ at which the convergence condition must be satisfied
(see Br\aa ting~\cite{Br} and Katz \& Katz \cite{KK11b}).

Cauchy's approach is based on two assumptions which can be stated in
modern terminology as follows:
\begin{enumerate}
\item
when you have a closed expression for a function, then its values at
``variable quantities'' (such as~$x=\frac{1}{n}$) are calculated by
using the same closed expression as at real values;
\item
to evaluate a function at a variable quantity generated by a sequence,
one evaluates term-by-term.
\end{enumerate}
Cauchy's strengthened condition amounts to requiring the error
$r_n(x)=f(x)-f_n(x)$ to become infinitesimal:
\begin{equation}
\label{62}
\text{if } x \text{ infinitesimal then } r_n(x) \text{ infinitesimal},
\end{equation}
which in the case of~$x$ given by~$(1/n)$ translates into the
requirement that~$r_n(1/n)$ tends to zero.

An epsilontic interpretation (in the context of an Archimedean
continuum) of Cauchy's 1821 and 1853 texts is untenable, as it
necessitates a {\em pair\/} of variables as in~\eqref{32}, where
Cauchy only used a {\em single one\/}, namely~$x$, but one drawn from
a ``thicker'' continuum including infinitesimals.  Namely, Cauchy
draws the points to be evaluated at from an infinitesimal-enriched
continuum.

We will refer to an infinitesimal-enriched continuum as a {\em
Bernoullian continuum\/}, or a ``B-continuum'' for short, in an
allusion to Johann Bernouilli.%
\footnote{\label{ber1}Bernoulli was the first to use infinitesimals in
a systematic fashion as a foundational concept, Leibniz himself having
employed both a syncategorematic and a true infinitesimal approach.
The pair of approaches in Leibniz are discussed by Bos \cite[item 4.2,
p.~55]{Bos}; see also Appendix~\ref{rival2}, footnote~\ref{ber2}.}

A null sequence such as~$1/n$ ``becomes" an infinitesimal, in Cauchy's
terminology.  Evaluating at points of a Bernoullian continuum makes it
possible to express uniform convergence in terms of a single
variable~$x$ rather than a pair~$(x,n)$.

Once one acknowledges that there are {\em two\/} variables in the
traditional epsilontic definition of uniform continuity and uniform
convergence, it becomes untenable to argue that the condition Cauchy
introduced was epsilontic uniform convergence.  A historian who
describes Cauchy's condition as uniform convergence, must acknowledge
that the definition involves an infinitesimal-enriched continuum, at
variance with Boyer's interpretation.

In Appendix~\ref{rival2}, Subsection~\ref{a1} we present a parallel
distinction between continuity and uniform continuity, where a similar
distinction in terms of the number of variables is made.

\section{Did Weierstrass succeed in eliminating infinitesimals?}
\label{six}


Did Weierstrass succeed in eliminating infinitesimals from
mathematics?

\renewcommand{\arraystretch}{1.3}
\begin{table}
\[
\begin{tabular}[t]
{ | p{.4in} || p{1.45in} | p{2.7in} | } \hline years & author &
contribution \\ \hline\hline 1821 & Cauchy & Infinitesimal definition
of continuity \\ \hline 1827 & Cauchy & Infinitesimal delta function
\\ \hline 1829 & Cauchy & Defined ``order of infinitesimal'' in terms
of rates of growth of functions \\ \hline 1852 & Bj\"orling & Dealt
with convergence at points ``indefinitely close'' to the limit \\
\hline 1853 & Cauchy & Clarified hypothesis of ``sum theorem'' by
requiring convergence at infinitesimal points \\ \hline 1870-1900 &
Stolz, du~Bois-Reymond, Veronese, and others & Infinitesimal-enriched
number systems defined in terms of rates of growth of functions\\
\hline 1902 & Emile Borel & Elaboration of du Bois-Reymond's system \\
\hline 1910 & G.~H.~Hardy & Provided a firm foundation for du
Bois-Reymond's orders of infinity \\ \hline 1926 & Artin--Schreier &
Theory of real closed fields \\ \hline 1930 & Tarski & Existence of
ultrafilters \\ \hline 1934 & Skolem & Nonstandard model of arithmetic
\\ \hline 1948 & Edwin Hewitt & Ultrapower construction of hyperreals
\\ \hline 1955 & \Los{} & Proved \Los's theorem forshadowing the
transfer principle \\ \hline 1961, 1966 & Abraham Robinson &
Non-Standard Analysis \\ \hline 1977 & Edward Nelson & Internal Set
Theory \\ \hline
\end{tabular}
\]
\caption{\textsf{Timeline of modern infinitesimals from Cauchy to
Nelson}.}
\label{heuristic}
\end{table}
\renewcommand{\arraystretch}{1}

The persistent idea that infinitesimals have been ``eliminated'' by
the great triumvirate of Cantor, Dedekind, and Weierstrass was soundly
refuted by Ehrlich \cite{Eh06}.  Ehrlich documents a rich and
uninterrupted chain of work on non-Archimedean systems, or what we
would call a Bernoullian continuum.  Some key developments in this
chain are listed in Table~\ref{heuristic} (see \cite{BK} for more
details).

The elimination claim can only be understood as an oversimplification
by Weierstrass's followers, who wish to view the history of analysis
as a triumphant march toward the radiant future of Weierstrassian
epsilontics.  Such a view of history is rejected by H.~Putnam who
comments on the introduction of the methods of the differential and
integral calculus by Newton and Leibniz in the following terms:
\begin{quote}
If the epsilon-delta methods had not been discovered, then
infinitesimals would have been postulated entities (just as
`imaginary' numbers were for a long time).  Indeed, this approach to
the calculus--enlarging the real number system--is just as consistent
as the standard approach, as we know today from the work of Abraham
Robinson [\ldots] If the calculus had not been `justified' Weierstrass
style, it would have been `justified' anyway (Putnam \cite{Pu75}).
\end{quote}
In short, there is a cognitive bias inherent in a postulation in an
inevitable outcome in the evolution of a scientific discipline.

The study of cognitive bias has its historical roots in Francis
Bacon's proposed classification of what he called {\em idola\/} (a
Latin plural) of several kinds.  He described these as things which
obstructed the path of correct scientific reasoning.  Of particular
interest to us are his {\em Idola fori\/} (``Illusions of the
Marketplace'': due to confusions in the use of language and taking
some words in science to have meaning different from their common
usage); and {\em Idola theatri\/} (``Illusions of the Theater'': the
following of academic dogma and not asking questions about the world),
see Bacon \cite{Bac}.

Completeness, continuity, continuum, Dedekind ``gaps'': these are
terms whose common meaning is frequently conflated with their
technical meaning.  In our experience, explaining
infinitesimal-enriched extensions of the reals to an epsilontically
trained mathematician typically elicits a puzzled reaction on her
part: ``But aren't the real numbers already {\em complete\/} by virtue
of having filled in all the {\em gaps\/} already?''

This question presupposes an academic dogma, viz., that there is a
single coherent conception of the continuum, and it is a complete,
Archimedean ordered field.  This dogma has recently been challenged.
Numerous possible conceptions of the continuum range from
S.~Feferman's predicative conception of the continuum \cite{Fef09}, to
F.~William Lawvere's \cite{Law} and J.~Bell's conception in terms of
an intuitionistic topos \cite{Be08}, \cite{Be09}, \cite{Be09b}.

To illustrate the variety of possible conceptions of the continuum,
note that traditionally, mathematicians have considered at least two
different types of continua.  These are Archimedean continua, or
A-continua for short, and infinitesimal-enriched (Bernoulli) continua,
or B-continua for short.  Neither an A-continuum nor a B-continuum
corresponds to a unique mathematical structure (see
Table~\ref{continuity}).  Thus, we have two distinct implementations
of an A-continuum:
\begin{itemize}
\item
the real numbers (or Stevin numbers),%
\footnote{\label{f31}See Section~\ref{stev}.}
in the context of classical logic (incorporating the law of excluded
middle);
\item
Brouwer's continuum built from ``free-choice sequences'', in the
context of intuitionistic logic.
\end{itemize}

\renewcommand{\arraystretch}{1.3}
\begin{table}
\[
\begin{tabular}[t]
{ | p{.9in} || p{1in} | p{.9in} | p{.5in} | p{.75in} | p{.5in} |}
\hline & Archimedean & Bernoullian \\ \hline\hline classical &
Stevin's \mbox{continuum}$^{\tiny \ref{f31}}$ & Robinson's
continuum$^{\tiny \ref{f4bis}}$ \\ \hline intuitionistic & Brouwer's
continuum & Lawvere's continuum$^{\tiny \ref{f5}}$ \\ \hline
\end{tabular}
\]
\caption{\textsf{Varieties of continua, mapped out according to a pair
of binary parameters: classical/intuitionistic and
Archimedean/Bernoullian}.}
\label{continuity}
\end{table}
\renewcommand{\arraystretch}{1}

John L. Bell describes a distinction within the class of an
infinitesimal-enriched B-continuum, in the following terms.
Historically, there were two main approaches to such an enriched
continuum, one by Leibniz, and one by B.~Nieuwentijt, who favored
nilpotent (nilsquare) infinitesimals whose squares are zero.
Mancosu's discussion of Nieuwentijt in \cite[chapter 6]{Ma96} is the
only one to date to provide a contextual understanding of
Nieuwentijt's thought (see also Mancosu and Vailati \cite{MV}).
J.~Bell notes:
\begin{quote}
Leibnizian infinitesimals (differentials) are realized in
[A.~Robinson's] nonstandard analysis,%
\footnote{\label{f4bis}More precisely, the Hewitt-{\L}o{\'s}-Robinson
continuum; see Appendix~\ref{rival2}.}
and nilsquare infinitesimals in [Lawvere's] smooth infinitesimal
analysis~(Bell \cite{Be08, Be09}).
\end{quote}
The latter theory relies on intuitionistic logic.%
\footnote{\label{f5}Lawvere's infinitesimals rely on a
category-theoretic framework grounded in intuitionistic logic (see
J.~Bell~\cite{Be09}).}
An implementation of an infinitesimal-enriched continuum was developed
by P.~Giordano (see \cite{Gio10b, Gio10a}), combining elements of both
a classical and an intuitionistic continuum.  The Weirstrassian
continuum is but a single member of a diverse family of concepts.

\section{Did Dedekind discover the essence of continuity?}

\label{dede}


Did Dedekind discover ``the essence of continuity", and is such
essence captured by his cuts?

In Dedekind's approach, the ``essence'' of continuity amounts to the
numerical assertion that two non-rational numbers should be equal if
and only if they define the same Dedekind cut%
\footnote{We will ignore the slight technical complication arising
from the fact that there are two ways of defining the Dedekind cut
associated with a {\em rational\/} number.}
on the rationals.  Dedekind formulated his ``essence of continuity''
in the context of the geometric line in the following terms:
\begin{quote}
If all points of the straight line fall into two classes such that
every point of the first class lies to the left of every point of the
second class, then there exists one and only one point which produces
this division of all points into two classes, this severing of the
straight line into two portions (Dedekind \cite[p. 11]{Ded63}).
\end{quote}
We will refer to this essence as the \emph{geometric essence of
continuity}.%
\footnote{Note that the geometric essence of continuity necessarily
fails over an ordered non-Archimedean field $F$.  Thus, defining
infinitesimals as elements of $F$ violating the traditional
Archimedean property, we can start with the cut of $F$ into positive
and negative elements, and then modify this cut by assigning all
infinitesimals to, say, the negative side.  Such a cut does not
correspond to an element of $F$.}
Dedekind goes on to comment on the epistemological status of this
statement of the essence of continuity:
\begin{quote}
[\dots] I may say that I am glad if every one finds the above
principle so obvious and so in harmony with his own ideas of a line;
for I am utterly unable to adduce any proof of its correctness, nor
has any one the power (ibid.).
\end{quote}
Having enriched the domain of rationals by adding irrationals, numbers
defined completely by cuts not produced by a rational, Dedekind
observes:
\begin{quote}
From now on, therefore, to every definite cut there corresponds a
definite rational or irrational number, and we regard two numbers as
different or unequal always and only when they correspond to
essentially different cuts (Dedekind \cite[p.~15]{Ded63}).
\end{quote}
By now, Dedekind postulates that two numbers should be equal ``always
and only'' [i.e., if and only if] they define identical cuts on the
{\em rational\/} numbers.  Thus, Dedekind postulates that there should
be ``one and only one'' number which produces such a division.
Dedekind clearly presents this as an exact arithmetic analogue to the
geometric essence of continuity.  We will refer to such a postulate as
the \emph{rational essence of continuity}.

Dedekind's postulation of rational essence is not accompanied by
epistemological worries as was his geometric essence a few pages
earlier.  Yet, rational essence entails a suppression of
infinitesimals: a pair of distinct non-rational numbers can define the
same Dedekind cut on~$\Q$, such as~$\pi$ and~$\pi+h$ with~$h$
infinitesimal; but one cannot have such a pair if one postulates the
rational essence of continuity, as Dedekind does.

Dedekind's technical work on the foundations of analysis has been
justly celebrated (see D.~Fowler \cite{Fo}).  Whereas everyone before
Dedekind had {\em assumed\/} that operations such as powers, roots,
and logarithms can be performed, he was the first to show how these
operations can be defined, and shown to be coherent, in the realm of
the real numbers (see Dedekind \cite[\S 6]{Ded}).

Meanwhile, the nature of his interpretive speculations about what does
or does not constitute the ``essence" of continuity, is a separate
issue.  For over a hundred years now, many mathematicians have been
making the assumption that space conforms to Dedekind's idea of ``the
essence of continuity", which in arithmetic translates into the
numerical assertion that two numbers should be equal if and only if
they define the same Dedekind cut on the rationals.  Such an
assumption rules out infinitesimals.  In the context of the hyperreal
number system, it amounts to an application of the standard part
function (see Appendix~\ref{rival2}), which forces the collapse of the
entire halo (cluster of infinitely close, or adequal, points) to a
single point.

The formal/axiomatic transformation of mathematics accomplished at the
end of the 19th century created a specific foundational framework for
analysis.  Weierstrass's followers raised a philosophical prejudice
against infinitesimals to the status of an axiom.  Dedekind's
``essence of continuity'' was, in essence, a way of steamrolling
infinitesimals out of existence.

In 1977, E. Nelson~\cite{Ne} created a set-theoretic framework
(enriching ZFC) where the usual construction of the reals produces a
number system containing entities that behave like infinitesimals.
Thus, the elimination thereof was not the only way to achieve rigor in
analysis as advertized at the time, but rather a decision to develop
analysis in just one specific way.

\section{Who invented Dirac's delta function?}

A prevailing sentiment today is that one of the spectacular successes
of the rigorous analysis was the justification of delta functions,
originally introduced informally by to P.~Dirac (1902--1984), in terms
of distribution theory.  But was it originally introduced informally
by Dirac?

In fact, Fourier \cite{Fou} and Cauchy exploited the ``Dirac'' delta
function over a century earlier.  Cauchy defined such functions in
terms of infinitesimals (see L\"utzen \cite{Lut82} and
Laugwitz~\cite{Lau92}).  A function of the type generally attributed
to Dirac was specifically described by Cauchy in 1827 in terms of
infinitesimals.  More specifically, Cauchy uses a unit-impulse,
infinitely tall, infinitely narrow delta function, as an integral
kernel.  Thus, in 1827, Cauchy used infinitesimals in his definition
of a ``Dirac'' delta function \cite[p.~188]{Ca27}.  Here Cauchy uses
infinitesimals~$\alpha$ and~$\epsilon$, where~$\alpha$ is, in modern
terms, the ``scale parameter'' of the ``Cauchy distribution'',
whereas~$\epsilon$ gives the size of the interval of integration.
Cauchy wrote:
\begin{quote}
Moreover one finds, denoting by~$\alpha$,~$\epsilon$ two infinitely
small numbers,
\begin{equation}
\label{151}
\frac{1}{2} \int_{a-\epsilon}^{a+\epsilon} F(\mu) \frac{\alpha \;
d\mu}{\alpha^2 + (\mu-a)^2} = \frac{\pi}{2} F(a)
\end{equation}
\end{quote}
(see Cauchy \cite[p.~188]{Ca27}).  Such a formula extracts the value
of a function~$F$ at a point~$a$ by integrating~$F$ against a delta
function defined in terms of an infinitesimal parameter~$\alpha$ (see
and Laugwitz \cite[p.~230]{Lau89}).  The expression
\[
\frac{\alpha}{\alpha^2 + (\mu-a)^2}
\]
(for real~$\alpha$) is known as the {\em Cauchy distribution\/} in
probability theory.  The function is called the probability density
function.  A Cauchy distribution with an infinitesimal scale parameter
produces a function with Dirac-delta function behavior, exploited by
Cauchy already in 1827 in work on Fourier series and evaluation of
singular integrals.

\section{Is there continuity between Leibniz and Robinson?}


Is there continuity between historical infinitesimals and Robinson's
theory?

Historically, infinitesimals have often been represented by null
sequences.  Thus, Cauchy speaks of a variable quantity as becoming an
infinitesimal in 1821, and his variable quantities from that year are
generally understood to be sequences of discrete values (on the other
hand, in his 1823 he used continuous variable quantities).
Infinitesimal-enriched fields can in fact be obtained from sequences,
by means of an ultrapower construction, where a null sequence
generates an infinitesimal.  Such an approach was popularized by
Luxemburg \cite{Lu62} in 1962, and is based on the work by E. Hewitt
\cite{Hew} from 1948.  Even in Robinson's approach~\cite{Ro66} based
on the compactness theorem, a null sequence is present, though
well-hidden, namely in the countable collection of
axioms~$\epsilon<\frac{1}{n}$.  Thus, null sequences provide both a
cognitive and a technical link between historical infinitesimals
thought of as variable quantities taking discrete values, on the one
hand, and modern infinitesimals, on the other (see Katz \& Tall
\cite{KT}).

Leibniz's heuristic \emph{law of continuity} was implemented
mathematically as \Los's theorem and later as the \emph{transfer
principle} over the hyperreals (see Appendix~\ref{rival2}), while
Leibniz's heuristic law of homogeneity (see Leibniz
\cite[p.~380]{Le10b}) and Bos \cite[p.~33]{Bos}) was implemented
mathematically as the standard part function (see Katz and Sherry
\cite{KS}).

\section{Is Lakatos' take on Cauchy tainted by Kuhnian relativism?}

Does Lakatos's defense of infinitesimalist tradition rest upon an
ideological agenda of Kuhnian relativism and fallibilism, inapplicable
to mathematics?

G.~Schubring summarizes fallibilism as an
\begin{quote}
enthusiasm for revising traditional beliefs in the history of science
and reinterpreting the discipline from a theoretical, epistemological
perspective generated by Thomas Kuhn's (1962) work on the structure of
scientific revolutions.  Applying Popper's favorite keyword of
fallibilism, the statements of earlier scientists that historiography
had declared to be false were particularly attractive objects for such
an epistemologically guided revision (Schubring
\cite[p.~431--432]{Sch}).
\end{quote}
Schubring then takes on Lakatos in the following terms:
\begin{quote}
The philosopher Imre Lakatos (1922-1972)%
\footnote{The dates given by Schubring are incorrect.  The correct
dates are 1922-1974.}
was responsible for introducing these new approaches into the history
of mathematics.  One of the examples he analyzed and published in 1966
received a great deal of attention: Cauchy's theorem and the problem
of uniform convergence.  Lakatos refines Robinson's approach by
claiming that Cauchy's theorem had also been correct at the time,
because he had been working with infinitesimals (ibid.).
\end{quote}
However, Schubring's summary of the philosophical underpinnings of
Lakatos' interpretation of Cauchy's sum theorem is not followed up by
an analysis of Lakatos's position (see \cite{La76, La78}).  It is as
if Schubring felt that labels of ``Kuhnianism'' and ``fallibilism''
are sufficient grounds for dismissing a scholar.  Schubring proceeds
similarly to dismiss Laugwitz's reading of Cauchy as ``solipsistic"
\cite[p.~434]{Sch}.  Schubring accuses Laugwitz of interpreting
Cauchy's conceptions as
\begin{quote}
some hermetic closure of a {\em private\/} mathematics (Schubring
\cite[p.~435]{Sch}) [emphasis in the original--the authors];
\end{quote}
as well as being ``highly anomalous or isolated'' \cite[p.~441]{Sch}.

The fact is that Laugwitz is interpreting Cauchy's words according to
their plain meaning (see \cite{Lau87, Lau97}), as revealed by looking,
as Kuhn would suggest, at the context in which they occur.  The
context strongly recommends taking Cauchy's infinitesimals at face
value, rather than treating them as a sop to the management.  The
burden of proof falls upon Schubring to explain why the triumvirate
interpretation of Cauchy is not ``solipsistic'', ``hermetic'', or
``anomalous''.  The latter three modifiers could be more
apppropriately applied to Schubring's own interpretation of Cauchy's
infinitesimals as allegedly involving a {\em compromise\/} with rigor,
allegedly due to {\em tensions\/} with the management of the {\em
Ecole polytechnique\/}.  Schubring's interpretation is based on
Cauchy's use of the term {\em concilier\/} in Cauchy's comment on the
first page of his {\em Avertissement\/}:
\begin{quote}
Mon but principal a \'et\'e de concilier la rigueur, dont je m'\'etais
fait une lois dans mon {\em Cours d'Analyse\/}, avec la simplicit\'e
qui r\'esulte de la consid\'eration directe des quantit\'es infiniment
petites (Cauchy \cite[p.~10]{Ca23}).
\end{quote}
Let us examine Schubring's logic of conciliation.  A careful reading
of Cauchy's {\em Avertissement\/} in its entirety reveals that Cauchy
is referring to an altogether different source of tension, namely his
rejection of some of the procedures in Lagrange's {\em M\'ecanique
analytique\/} \cite{Lag} as unrigorous, such as Lagrange's principle
of the ``generality of algebra''.  While rejecting the ``generality of
algebra'' and Lagrange's flawed method of power series, Cauchy was
able, as it were, to sift the chaff from the grain, and retain the
infinitesimals endorsed in the 1811 edition of the {\em M\'ecanique
analytique\/}.  Indeed, Lagrange opens his treatise with an
unequivocal endorsement of infinitesimals.  Referring to the system of
infinitesimal calculus, Lagrange writes:
\begin{quote}
Lorsqu'on a bien con\c cu l'esprit de ce syst\`eme, et qu'on s'est
convaincu de l'exactitude de ses r\'esultats par la m\'ethode
g\'eom\'etrique des premi\`eres et derni\`eres raisons, ou par la
m\'ethode analytique des fonctions d\'eriv\'ees, on peut employer les
infiniment petits comme un instrument s\^ur et commode pour abr\'eger
et simplifier les d\'emonstrations%
\footnote{``Once one has duly captured the spirit of this system
[i.e., infinitesimal calculus], and has convinced oneself of the
correctness of its results by means of the geometric method of the
prime and ultimate ratios, or by means of the analytic method of
derivatives, one can then exploit the infinitely small as a reliable
and convenient tool so as to shorten and simplify proofs''
(Lagrange).}
(Lagrange \cite[p.~iv]{Lag}).
\end{quote}

Lagrange describes infinitesimals as dear to a scientist, being
reliable and convenient.  In his {\em Avertissement\/}, Cauchy retains
the infinitesimals that were also dear to Lagrange, while criticizing
Lagrange's ``generality of algebra'' (see~\cite{KK11a} for details).

It's useful here to evoke the use of the term ``concilier" by Cauchy's
teacher Lacroix.  Gilain quotes Lacroix in 1797 to the effect that
\begin{quote}
``lorsqu'on veut concilier la rapidit\'e de l'exposition avec
l'exactitude dans le langage, la clart\'e dans les principes, [\ldots],
je pense qu'il convient d'employer la m\'ethode des limites"
(p.XXIV).''  \cite[footnote~20]{Gi}.
\end{quote}
Here Lacroix, like Cauchy, employs ``concilier'', but in the context
of discussing the {\em limit\/} notion.  Would Schubring's logic of
conciliation dictate that Lacroix developed a compromise notion of
limit, similarly with the sole purpose of accomodating the management
of the {\em Ecole\/}?

Why are Lakatos and Laugwitz demonized, rather than analyzed, by
Schubring?  We suggest that the act of contemplating for a moment the
idea that Cauchy's infinitesimals can be taken at face value is
unthinkable to a triumvirate historian, as it would undermine the
epsilontic Cauchy-Weierstrass tale that the received historiography is
erected upon.  The failure to appreciate the Robinson-Lakatos-Laugwitz
interpretation, according to which infinitesimals are mainstream
analysis from Cauchy onwards, is symptomatic of a narrow
Archimedean-continuum vision.

\appendix

\section{Rival continua}
\label{rival2}

This appendix summarizes a 20th century implementation of an
alternative to an Archimedean continuum, namely an
infinitesimal-enriched continuum.  Such a continuum is not to be
confused with incipient notions of such a continuum found in earlier
centuries in the work of Fermat, Leibniz, Euler, Cauchy, and others.

Johann Bernoulli was one of the first to exploit infinitesimals in a
systematic fashion as a foundational tool in the calculus.%
\footnote{\label{ber2}See footnote~\ref{ber1} for a comparison with
Leibniz.}
We will therefore refer to such a continuum as a Bernoullian
continuum, or B-continuum for short.

\subsection{Constructing the hyperreals}

Let us start with some basic facts and definitions.  Let
\[
(\mathbb{R},+,\cdot,0,1,<)
\]
be the field of real numbers, let~$\mathcal{F}$ be a fixed
nonprincipal ultrafilter on~$\mathbb{N}$ (the existence of such was
established by Tarski \cite{Tar}).  The relation~$\equiv$
\mbox{defined by}
\[
(r_n)\equiv{(s_n)}\leftrightarrow_{\text{def}}
\{n\in{\mathbb{N}}:r_n=s_n\}\in{\mathcal{F}}
\]
is an equivalence relation on the set~$\mathbb{R}^\mathbb{N}$.  The
set of hyperreals~$\RRR$, or the B-continuum for short, is the
quotient set
\[
\RRR=_{\text{def}}\mathbb{R}^\mathbb{N}/_{\equiv}.
\]
Addition, multiplication and order of hyperreals are defined by
\[[
(r_n)] + [(s_n)]=_{\text{def}}[(r_n+s_n)],\ \ \ \ \ [(r_n)] \cdot
[(s_n)]=_{\text{def}}[(r_n\cdot{s_n})],
\]
\[
[(r_n)]\prec{[(s_n)]}\leftrightarrow_{\text{def}}
{\{n\in{\mathbb{N}}:r_n<s_n\}\in{\mathcal{F}}}.
\]
The standard real number~$r$ is identified with equivalence class
$r^*$ of the constant sequence~$(r,r,\ldots\ )$, i.e.
$r^*=_{\text{def}}[(r,r,\ldots\ )]$.

The set~$\NNN$ of hypernaturals (mentioned in Section~\ref{four})%
\footnote{\label{hyper2}See footnote~\ref{hyper1}.}
is the subset of~$\RRR$ defined by
\[
[(r_n)]\in{\NNN}\leftrightarrow_{\text{def}}
\{n\in{\N}:r_n\in{\N}\}\in{\mathcal{F}}.
\]
In particular, each sequence of natural numbers~$(n_j)$ represents a
hypernatural number, i.e.~$[(n_j)]\in{\NNN}$.

The set of hypernaturals can be represented as a disjoint union
\[
\NNN= \{n^*:n\in{\N}\}\cup \NNN_{\infty},
\] 
where the set~$\{n^*:n\in{\N}\}$ is just a copy of the usual natural
numbers, and~$\NNN_{\infty}$ consists of infinite (sometimes called
``unlimited'') hypernaturals.  Each element of~$\NNN_{\infty}$ is
greater than every usual natural number, i.e.
\[
(\forall{H\in{\NNN_{\infty}}}) (\forall{n\in{\N}}) \, [H\succ n^*].
\]

\begin{theorem}
$(\mathbb{R}^*, +, \cdot, 0^*,1^*,\prec)$ is a non-Archimedean, real
closed field.
\end{theorem}

The set of infinitesimal hyperreals~$\Omega$ is defined by
\[
x\in{\Omega}\text{ if and only if } (\forall{\theta\in{\R_{+}}}) \,
[\,|x|\prec\theta^*\ ],
\]
where~$|x|$ stands for the absolute value of~$x$, which is defined as
in any ordered field.  We say that~$x$ is infinitely close to~$y$, and
write~$x\approx{y}$, if and only if~$x-y\in{\Omega}$.

To give some examples, the sequence~$(\tfrac{1}{n})$ represents a
positive infinitesimal~$[(\tfrac{1}{n})]$.  Next, let~$(r_n)\in\R^\N$
be a sequence of reals such that
$\lim\limits_{n\rightarrow\infty}r_n=0$, then~$(r_n)$ represents an
infinitesimal,%
\footnote{In this construction, every null sequence defines an
infinitesimal, but the converse is not necessarily true.  Modulo
suitable foundational material, one can ensure that every
infinitesimal is represented by a null sequence; an appropriate
ultrafilter (called a {\em P-point\/}) will exist if one assumes the
continuum hypothesis, or even the weaker Martin's axiom.  See Cutland
{\em et al\/} \cite{CKKR} for details.}
i.e. ~$[(r_n)]\in{\Omega}$.  And finally, sequence
$\left(\tfrac{(-1)^n}{n}\right)$
represents a nonzero
infinitesimal~$\left[\left(\tfrac{(-1)^n}{n}\right)\right]$, whose
sign depends on whether or not the set~$2\N$ is a member of the
ultrafilter.

The set of limited hyperreals~$\RRR_{<\infty}$ is defined by
\[
x\in{\mathbb{\RRR_{<\infty}}}\leftrightarrow_{\text{def}}{\exists{\theta\in
{\R_{+}}}}[\ |x|\prec\theta^*\ ],
\]
so that we have a disjoint union
\begin{equation}
\RRR= \RRR_{<\infty} \cup \RRR_{\infty},
\end{equation}
where~$\RRR_{\infty}$ consists of unlimited hyperreals (i.e., inverses
of nonzero infinitesimals).

\begin{theorem}[Standard Part Theorem]
\[
(\forall{x\in\text{\rm \RRR}_{<\infty}}) (\exists!{r\in{\R}}) \, [\,
r^*\approx{x}\ ].
\]
\end{theorem}

The unique real~$r$ such that~$\ r^*\approx{x}$ is called the standard
part of~$x$, and we write~$\text{st}(x)=r$.

Note that if a sequence $(r_n:n\in\N)$ happens to be Cauchy, one can
relate standard part and limit as follows:%
\footnote{This theme is pursued further by Giordano et
al. \cite{GK11}.}
\begin{equation}
\st([(r_n)])=\lim_{n\to\infty} r_n.
\end{equation}

\begin{theorem}
$\Omega$ is a maximal ideal of the ring~$(\text{\rm
\RRR}_{<\infty},+,\cdot,0^*,1^*)$, and the quotient ring is isomorphic
to the field of standard real numbers~$(\mathbb{R},+,\cdot,0,1)$.
\end{theorem}

Since the map
\[
\R\ni r\mapsto r^*\in{\RRR}
\]
is an order preserving morphism, we can treat the field of hyperreals
as an extension of standard reals and use the usual notation
$(\RRR,+,\cdot,0,1,<)$. Now, the map ``st'' sends each finite
point~$x\in \RRR_{<\infty}$, to the real point st$(x)\in \R$
infinitely close to~$x$, as follows:%
\footnote{This is the Fermat-Robinson standard part whose seeds are
found in Fermat's adequality, as well as in Leibniz's treanscendental
principle of homogeneity.}
\begin{equation*}
\xymatrix{\quad \RRR_{{<\infty}}^{~} \ar[d]^{{\rm st}} \\ \R}
\end{equation*}
Robinson's answer to Berkeley's {\em logical criticism\/} (see
D.~Sherry \cite{She87}) is to define the basic concept of the calculus
as
\begin{equation*}
\hbox{st} \left( \frac{\Delta y}{\Delta x} \right),
\end{equation*}
rather than the differential ratio~$\Delta y/\Delta x$ itself, as in
Leibniz.  Robinson comments as follows: ``However, this is a small
price to pay for the removal of an inconsistency'' (Robinson
\cite[p~266]{Ro66}).%
\footnote{However, as argued in \cite{KS}, an alleged inconsistency
was not there in the first place.  Briefly, Leibniz is able to employ
his transcendental law of homogeneity to the same effect as Robinson's
standard part function (see Bos 1974, \cite[p.~33]{Bos}).}

A sequence~$(r_n:n\in{\N})$ of real numbers can be extended to a
hypersequence~$(r_K:K\in\text{\rm \NNN})$ of hyperreals, indexed by
\emph{all} the hypernaturals, by setting
\[
r_K=[(r_{k_1},r_{k_2}, \ldots,r_{k_j},\ldots)],\ \ \ \mbox{where}\
K=[(k_1,k_2,\ldots,k_j,\ldots)].
\]

\begin{theorem}
Let~$(r_n)$ be a sequence of real numbers, and let~$ L\in{\R}$. Then
\[
\lim\limits_{n\rightarrow \infty}r_n= L \leftrightarrow
(\forall{H\in{\text{\rm \NNN}_{\infty}}}) \, [\,r_H\approx L^*].
\]
\end{theorem}

\subsection{Uniform continuity}
\label{a1}

We present a discussion of uniform continuity so as to supplement and
clarify the discussion of uniform {\em convergence\/} in
Section~\ref{sum}.  The idea of ``one variable {\em versus\/} two
variables" is a little easier to explain in the context of uniform
continuity.

Thus, the traditional definition of ordinary continuity on an interval
can be expressed in terms of a single variable~$x$ running through the
domain~$D_f$ of the function~$f$: namely,
\begin{quote}
for each~$x\in D_f$, \quad~$\lim f(x') = f(x)$ as~$x'$ tends to~$x$.
\end{quote}
Meanwhile, uniform continuity cannot be expressed in a similar way in
the traditional framework.  Namely, one needs a {\em pair\/} of
variables to run through the domain of~$f$:
\begin{quote}
$\forall x, y \in D_f$, \quad if~$|x-y|<\delta$ then
$|f(x)-f(y)|<\epsilon$
\end{quote}
(of course, this has to be prefaced by the traditional epsilon-delta
yoga).  Now the crucial observation is that in the context of a
B-continuum, one no longer needs a \emph{pair} of variables to define
uniform continuity.  Namely, it can be defined using a single
variable, by exploiting the notion of {\em microcontinuity\/} at a
point (see Gordon et al. \cite{GKK}).  We will use Leibniz's
symbol~$\adequal$ for the relation of being infinitely close.
Thus,~$f$ is called microcontinuous at~$x$ if
\begin{quote}
whenever~$y\adequal x$, also~$f(y) \adequal f(x)$.
\end{quote}
In terms of this notion, the uniform continuity of a real function~$f$
is defined in terms of its natural extension~$f^*$ to the hyperreals
as follows:
\begin{quote}
for all~$x \in D_{f^*}, \quad f^*$ is microcontinuous at~$x$.
\end{quote}
This sounds startlingly similar to the definition of continuity
itself, but the point is that microcontinuity is now required at every
point of the B-continuum, i.e., in the domain of~$f^*$ which is the
natural extension of the (real) domain of~$f$.

To give an example, the function~$f(x)=x^2$ fails to be uniformly
continuous on~$\R$ because of the failure of microcontinuity of its
natural extension~$f^*$ at a single infinite hyperreal~$H$.  The
failure of microcontinuity at~$H$ is checked as follows.  Consider the
infinitesimal~$e=\frac{1}{H}$, and the point~$H+e$ infinitely close
to~$H$.  To show that~$f^*$ is not microcontinuous at~$H$, we
calculate
\[
f^*(H+e)=(H+e)^2 = H^2 + 2He + e^2 = H^2 + 2 +e ^2 \adequal H^2 + 2.
\]
This value is not infinitely close to~$f^*(H)=H^2$, hence
microcontinuity fails at~$H$.  Thus the squaring function~$f(x)=x^2$
is not uniformly continuous on~$\R$.

We introduced the term ``microcontinuity'' (cf.~Gordon et al
\cite{GKK}) since there are two definitions of continuity, one using
infinitesimals, and one using epsilons.  The former is what we refer
to as microcontinuity.  It is given a special name to distinguish it
from the traditional definition of continuity.  Note that
microcontinuity at a non-standard hyperreal does not correspond to any
notion available in the epsilontic framework.  To give another
example, if we consider the function~$f$ given by $f(x)=\frac{1}{x}$
on the open interval~$(0,1)$ as well as its natural extension~$f^*$,
then~$f^*$ fails to be microcontinuous at a positive
infinitesimal~$e>0$.  It follows that~$f$ is not uniformly continuous
on~$(0,1)$.

\subsection{Pedagogical advantage of microcontinuity}

The expressibility of uniform continuity in terms of a condition on a
\emph{single} variable, as explained above, is a pedagogical advantage
of the microcontinuous approach.  The pedagogical difficulty of the
traditional two-variable definition in an epsilontic framework is
compounded by its multiple alternations of quantifiers, while the
hyperreal approach reduces the logical complexity of the definition by
two quantifiers (see, e.g., Keisler~\cite{Ke08b}).  The natural
hyperreal extension $f^*$ of a real function $f$ is, of course,
necessarily continuous at non-standard points, as well, by the
transfer principle; on the other hand, this type of continuity at a
non-standard point is of merely theoretical relevance in a calculus
classroom.  The relevant notion is that of microcontinuity, which
allows one to distinguish between the classical notions of continuity
and uniform continuity in a lucid way available only in the hyperreal
framework.  Similarly, the failure of uniform continuity can be
checked by a ``covariant''%
\footnote{See discussion at the end of Subsection~\ref{lutz}.}
(direct) calculation at a single non-standard point, whereas the
argument in an epsilontic framework is a bit of a ``contravariant''
multiple-quantifier tongue twister.

We are therefore puzzled by Hrbacek's dubious laments \cite{Hr05,
Hr10} of alleged ``pedagogical difficulties'' related to behavior at
non-standard points, directed at the framework developed by Robinson
and Keisler.  On the contrary, such a framework bestows a distinct
pedagogical advantage.%
\footnote{\label{connes}Another critic of Robinson's framework is
A.~Connes.  He criticizes Robinson's infinitesimals for being
dependent on non-constructive foundational material.  He further
claims it to be a {\em weakness\/} of Robinson's infinitesimals that
the results of calculations that employ them, do not depend on the
choice of the infinitesimal.  Yet, Connes himself develops a theory of
infinitesimals bearing a similarity to the ultrapower construction of
the hyperreals in that it also relies on sequences (more precisely,
spectra of compact operators).  Furthermore, he freely relies on such
results as the existence of the Dixmier trace, and the Hahn-Banach
theorem.  The latter results rely on similarly nonconstructive
foundational material.  Connes claims the independence of the choice
of Dixmier trace to be a {\em strength\/} of his theory of
infinitesimals in \cite[p.~6213]{Co95}.  Thus, both of Connes'
criticisms apply to his own theory of infinitesimals.  The
mathematical novelty of Connes' theory of infinitesimals resides in
the exploitation of Dixmier's trace, relying as it does on
non-constructive foundational material, thus of similar foundational
status to, for instance, the ultrapower construction of a
non-Archimedean extension of the reals (see also \cite{KK11d}).}

\subsection{Historical remarks}

Both the term ``hyper-real field'', and an ultrapower construction
thereof, are due to E.~Hewitt in 1948 (see \cite[p.~74]{Hew}).  In
1966, Robinson referred to the
\begin{quote}
theory of hyperreal fields (Hewitt [1948]) which \ldots can serve as
non-standard models of analysis \cite[p.~278]{Ro66}.
\end{quote}
The {\em transfer principle\/} is a mathematical implementation of
Leibniz's heuristic {\em law of continuity\/}: ``what succeeds for the
finite numbers succeeds also for the infinite numbers and vice versa''
(see Robinson~\cite[p.~262, 266]{Ro66} citing Leibniz 1701,
\cite{Le01d}).  The transfer principle, allowing an extension of every
first-order real statement to the hyperreals, is a consequence of the
theorem of J.~{\L}o{\'s} in 1955, see~\cite{Lo}, and can therefore be
referred to as a Leibniz-{\L}o{\'s} transfer principle.  A
Hewitt-{\L}o{\'s} framework allows one to work in a B-continuum
satisfying the transfer principle.

\subsection{Comparison with Cantor's construction}
To indicate some similarities between the ultrapower construction and
the so called Cantor's construction of the real numbers, let us start
with the field of rational numbers~$\Q$.  Let~$\Q^\N$ be the ring of
sequences of rational numbers.  Denote by
\begin{equation*}
\left( \Q^\N \right)_C
\end{equation*}
the subspace consisting of Cauchy sequences, and
let~$\mathcal{F}_{\!n\!u\!l\!l}\subset \left( \Q^\N \right)_C$ be the
subspace of all null sequences.  The reals are by definition the
quotient field
\begin{equation}
\label{real}
\R= \left. \left( \Q^\N \right)_C \right/ \mathcal{F}_{\!n\!u\!l\!l}.
\end{equation}
Meanwhile, an infinitesimal-enriched field extension of~$\Q$ may be
obtained by forming the quotient
\begin{equation*}
\left.  \Q^\N \right/ \mathcal{F}.
\end{equation*}
Here~$[(q_n )]$ maps to zero in the quotient if and only if one has
\[
\{ n \in \N : q_n = 0 \} \in {\mathcal{F}},
\]
where~$\mathcal{F}$, as above, is a fixed nonprincipal ultrafilter
on~$\N$.

\begin{figure}
\begin{equation*}
\xymatrix{ && \left( \left. \Q^\N \right/ \mathcal{F}
\right)_{<\infty} \ar@{^{(}->} [rr]^{} \ar@{->>}[d]^{\rm st} &&
\RRR_{<\infty} \ar@{->>}[d]^{\rm st} \\ \Q \ar[rr] \ar@{^{(}->} [urr]
&& \R \ar[rr]^{\simeq} && \R }
\end{equation*}
\caption{\textsf{An intermediate field~$\left. \Q^\N \right/
\mathcal{F}$ is built directly out of~$\Q$}}
\label{helpful}
\end{figure}

To obtain a full hyperreal field, we replace~$\Q$ by~$\R$ in the
construction, and form a similar quotient
\begin{equation}
\label{bee}
\RRR= \left.  \R^\N \right/ \mathcal{F}.
\end{equation}
We wish to emphasize the analogy with formula~\eqref{real} defining
the A-continuum.  We can treat both~$\R$ and~$\Q^\N / \mathcal{F}$ as
subsets of~$\RRR$.  Note that, while the leftmost vertical arrow in
Figure~\ref{helpful} is surjective, we have
\begin{equation*}
\left( \Q^\N / \mathcal{F} \right) \cap \R = \Q.
\end{equation*}

\subsection{Applications}
A more detailed discussion of the ultrapower construction can be found
in M.~Davis~\cite{Dav77} and Gordon, Kusraev, \&
Kutateladze~\cite{GKK}.  See also B\l aszczyk \cite{Bl} for some
philosophical implications.  More advanced properties of the
hyperreals such as saturation were proved later (see Keisler
\cite{Ke94} for a historical outline).  A helpful ``semicolon''
notation for presenting an extended decimal expansion of a hyperreal
was described by A.~H.~Lightstone~\cite{Li}.  See also P.~Roquette
\cite{Roq} for infinitesimal reminiscences.  A discussion of
infinitesimal optics is in K.~Stroyan \cite{Str},
J.~Keisler~\cite{Ke}, D.~Tall~\cite{Ta80}, and L.~Magnani and
R.~Dossena~\cite{MD, DM}.

Edward Nelson \cite{Ne} in 1977 proposed an axiomatic theory parallel
to Robinson's theory. A related theory was proposed by Hrb\'a\v{cek}
\cite{Hr} (who submitted a few months earlier and published a few
months later than Nelson).  Another axiomatic approach was proposed by
Benci and Di~Nasso \cite{BD}.  As Ehrlich \cite[Theorem~20]{Eh12}
showed, the ordered field underlying a maximal (i.e.,
\emph{On}-saturated) hyperreal field is isomorphic to J. H. Conway's
ordered field No, an ordered field Ehrlich describes as the
\emph{absolute arithmetic continuum}.

Infinitesimals can be constructed out of integers (see Borovik, Jin,
and Katz \cite{BJK}).  They can also be constructed by refining
Cantor's equivalence relation among Cauchy sequences (see Giordano \&
Katz \cite{GK11}).  A recent book by Terence Tao contains a discussion
of the hyperreals \cite[p. 209-229]{Tao}.

The use of the B-continuum as an aid in teaching calculus has been
examined by Tall \cite{Ta91}, \cite{Ta09a}; Ely \cite{El}; Katz and
Tall \cite{KT} (see also \cite{KK1, KK2}).  These texts deal with a
``naturally occurring'', or ``heuristic'', infinitesimal entity
$1-\text{`}0.999\ldots$' and its role in calculus pedagogy.%
\footnote{\label{per2}See footnote~\ref{per1} for Peirce's take on
$1-\text{`}0.999\ldots$'.}
Applications of the B-continuum range from the Bolzmann equation (see
L.~Arkeryd~\cite{Ar81, Ar05}); to modeling of timed systems in
computer science (see H.~Rust \cite{Rust}); Brownian motion, economics
(see R.~Anderson \cite{An76}); mathematical physics (see Albeverio
{\em et al.\/} \cite{Alb}); etc.

\section*{Acknowledgments}

We are grateful to M.~Barany, L.~Corry, N.~Guicciardini, and V.~Katz
for helpful comments.

{\bf Piotr B\l aszczyk} is Professor at the Institute of Mathematics,
Pedagogical University (Cracow, Poland).  He obtained degrees in
mathematics (1986) and philosophy (1994) from Jagiellonian University
(Cracow, Poland), and a PhD in ontology (2002) from Jagiellonian
University.  He authored \emph{Philosophical Analysis of Richard
Dedekind's memoir} Stetigkeit und irrationale Zahlen (2008,
Habilitationsschrift).  His research interest is in the idea of
continuum and continuity from Euclid to modern times.

\medskip

{\bf Mikhail G. Katz} is Professor of Mathematics at Bar Ilan
University, Ramat Gan, Israel.  Two of his joint studies with Karin
Katz were published in {\em Foundations of Science\/}: ``A Burgessian
critique of nominalistic tendencies in contemporary mathematics and
its historiography" and ``Stevin numbers and reality", online
respectively at

\url{http://dx.doi.org/10.1007/s10699-011-9223-1} and at

\url{http://dx.doi.org/10.1007/s10699-011-9228-9}

A joint study with Karin Katz entitled ``Meaning in classical
mathematics: is it at odds with Intuitionism?" {\em Intellectica\/}
\textbf{56} (2011), no.~2, 223-302 may be found at
\url{http://arxiv.org/abs/1110.5456}

A joint study with A. Borovik and R. Jin entitled ``An integer
construction of infinitesimals: Toward a theory of Eudoxus
hyperreals'' is due to appear in \emph{Notre Dame Journal of Formal
Logic} \textbf{53} (2012), no.~4.

A joint study with David Sherry entitled ``Leibniz's infinitesimals:
Their fictionality, their modern implementations, and their foes from
Berkeley to Russell and beyond'' is due to appear in
\emph{Erkenntnis}.

A joint study with David Tall, entitled ``The tension between
intuitive infinitesimals and formal mathematical analysis", appeared
as a chapter in a book edited by Bharath Sriraman, see

\noindent
\url{http://www.infoagepub.com/products/Crossroads-in-the-History-of-Mathematics}

\medskip

{\bf David Sherry} is fortunate to be professor of philosophy at
Northern Arizona University.

\end{document}